\documentclass[preprint,12pt]{elsarticle}

\usepackage{amssymb}
\usepackage{amscd,amsmath,amsfonts,graphicx}
\usepackage{amsthm}

\begin{document}

\baselineskip=22pt \centerline{\Large \bf Recent Developments on the
Moment Problem} \vspace{1cm} \centerline{Gwo Dong Lin}
\centerline{Academia Sinica, Taiwan} \vspace{0.7cm}\noindent {\bf
Abstract.} We consider univariate distributions with finite moments
of all positive orders. The moment problem is to determine whether
or not a given distribution is uniquely determined by the sequence
of its moments. There is a huge literature on this
classical topic. In this survey, we will  focus only on the recent
developments on the {\it checkable} moment-(in)determinacy criteria
including Cram\'er's condition, Carleman's condition, Hardy's
condition, Krein's condition and the growth rate of moments, which
help us solve the problem more easily. Both Hamburger and Stieltjes
cases are investigated. The former is concerned with distributions
on the whole real line, while the latter deals only with
distributions on the right half-line. Some new results or new simple
(direct) proofs of previous criteria are provided. Finally, we
review the most recent moment problem for products of independent
random variables with different distributions, which occur naturally
in stochastic modelling of complex random phenomena.

\vspace{0.7cm} \hrule \medskip \vspace{0.2cm}\noindent{\bf 2010 AMS
Mathematics Subject Classifications}:  60E05, 44A60.
{\bf Postal address:} Gwo Dong Lin, Institute of Statistical
Science, Academia Sinica, Taipei
11529, Taiwan, R.O.C. (E-mail: gdlin@stat.sinica.edu.tw)\\

\newpage \noindent {\bf 1. Introduction}

The moment problem is a classical topic over one century old (Stieltjes
1894/1895, Kjeldsen 1993, Fischer 2011, pp.\,157--168). We start
with the definition of the moment determinacy of distributions.
 Let $X$ be a random variable with distribution $F$ (denoted $X\sim F$)
 and have finite moments
 $m_k={\bf E}[X^k]$ for all  $k=1,2,\ldots;$
 namely, the absolute moment $\mu_k={\bf E}[|X|^k]<\infty$ for all
 positive integers $k$.
If $F$ is {uniquely determined} by the sequence of its moments
$\{m_k\}_{k=1}^{\infty}$, we say that $F$ is {moment-determinate}
(in short, $F$ is {M-det}, or $X$ is M-det); otherwise, we say that
$F$  is moment-indeterminate ($F$ is {M-indet}, or $X$ is M-indet).

 The moment problem is to determine whether
or not a given distribution $F$ is M-det. Roughly speaking, there
are two kinds of moment problems:  Stieltjes (1894/1895) moment
problem deals with nonnegative random variables only, while
Hamburger (1920/1921) moment problem treats all random variables
taking values in the whole real line.

   We recall first two important facts:\\
{\bf Fact A.} It is possible that a nonnegative random variable $X$
is M-det in the Stieltjes sense, but M-indet in the Hamburger sense
(Akhiezer 1965, p.\,240). This happens only for some  {\it discrete}
nonnegative random variables
with a positive mass at zero (Chihara 1968). \\
{\bf Fact B.} If a distribution $F$ is {M-indet}, then there are
infinitely many (i) absolutely continuous distributions, (ii) purely
discrete distributions and (iii) singular continuous distributions
all having the same moment sequence as $F$ (Berg 1998, Berg and
Christensen 1981).

  One good reason to study the moment problem was given in
 Fr$\acute{\hbox{e}}$chet and Shohat's (1931) Theorem stated below.  {Simply speaking,} for a
given
 sequence of random variables $X_n\sim
F_n,\ n=1,2,\ldots,$ with finite moments $m_k^{(n)}={\bf E}[X_n^k]$
for all positive integers $k$, the moment convergence
($\lim_{n\rightarrow\infty} m_k^{(n)}=m_k\ \forall k$) does not
guarantee the weak convergence of distributions
$\{F_n\}_{n=1}^{\infty}$ ($F_n\stackrel{\scriptsize
\hbox{w}}{\rightarrow} F$ as $n\to \infty$) unless the limiting
distribution $F$ is M-det. Therefore, the M-(in)det property is one
of the important fundamental properties we have to know about a
given distribution.

\noindent {\bf Fr$\acute{\hbox{e}}$chet and Shohat's (1931)
Theorem.} Let the distribution functions $F_n$ possess finite
moments $m_k^{(n)}$ for $k=1,2,\ldots$ and $n=1,2,\ldots.$ Assume
further that the limit $m_k=\lim_{n\rightarrow\infty} m_k^{(n)}$
exists (and is finite) for
each $k$. Then\\
(i) the limits $\{m_k\}_{k=1}^{\infty}$ are the moment sequence of a distribution function, say $F$;\\
(ii) if the limit $F$ given by (i) is {M-det}, $F_n$ converges to
$F$ weakly as $n\to\infty$.

Necessary and sufficient conditions for the M-det property of
distributions exist in the literature (see, e.g., Akhiezer 1961,
Shohat and Tamarkin 1943, and Berg et al. 2002), but these
conditions are not easily checkable in general. In this survey, we
will focus only on the {\it checkable} M-(in)det criteria for
distributions rather than the collection of all specific examples.

In Sections 2 and 3, we review respectively the moment determinacy
and moment indeterminacy criteria including Cram\'er's condition,
Carleman's condition, Hardy's condition, Krein's condition and the
growth rate of moments. Some criteria are old, but others are
recent. New (direct) proofs for some criteria are  provided. To
amend some previous proofs in the literature,  two lemmas (Lemmas 3
and 4) are given for the first time. We consider in Section 4 the
recently formulated Stieltjes classes for M-indet absolutely
continuous distributions. Section 5 is devoted to the converses to
the previous M-(in)det criteria for distributions. Finally, in
Section 6 we review the most recent results about the moment problem
for products
of independent random variables with different distributions.\medskip\\
\noindent {\bf 2. Checkable Criteria for Moment Determinacy}

 In this section we consider the checkable criteria for moment determinacy
 of random variables or distributions. We treat
first the Hamburger case because it is more popular than the
 Stieltjes case.  Let {$X\sim F$ on the whole real line ${\mathbb R}=(-\infty,\infty)$} with finite moments
$m_k={\bf E}[X^k]$  and absolute moment  $\mu_k={\bf E}[|X|^k]$
for all positive integers $k$. For convenience, we
define the following statements, in which `h' stands for `Hamburger'.\medskip\\
(h1) $\frac{m_{2(k+1)}}{m_{2k}}={\cal O}((k+1)^2)={\cal O}(k^2)$ as $k\to\infty$. \\
{(h2)} $X$ has a moment generating function (mgf), i.e., ${\bf
E}[e^{t{X}}]<\infty$ for all $t\in (-c,c)$, where\\
\indent ~~$c>0$ is a constant (Cram\'er's condition);
equivalently, ${\bf E}[e^{t{|X|}}]<\infty$ for $0\le t<c$.\\
{(h3)}
$\limsup_{k\to\infty}\frac{1}{2k}m_{2k}^{1/(2k)}<\infty.$ \\
{(h4)} $\limsup_{k\to\infty}\frac{1}{k}\mu_{k}^{1/k}<\infty.$\\
{(h5)} $m_{{2k}}={\cal O}((2k)^{2k})$ as
$k\to\infty.$\\
{(h6)} $m_{{2k}}\le c_0^k\,(2k)!,\ k=1,2,\ldots,$ for some
constant
$c_0>0.$ \\
(h7) ${C}[F]\equiv\sum_{k=1}^{\infty}m_{2k}^{-1/(2k)}=\infty$
(Carleman's (1926) condition).\\
(h8) $X$ is M-det on ${\mathbb R}$.

\noindent {\bf Theorem 1.} Under the above settings, if $X\sim F$ on
${\mathbb R}$ satisfies one of the conditions (h1) through (h7),
then $X$ is M-det on ${\mathbb R}$. Moreover, (h1) implies (h2),
(h2) through (h6) are
equivalent, and (h6) implies (h7). In other words,  the following chain of implications holds:\\
\centerline{ (h1) $\Longrightarrow$ (h2) $\Longleftrightarrow$ (h3)
$\Longleftrightarrow$ (h4) $\Longleftrightarrow$ (h5)
$\Longleftrightarrow$ (h6) $\Longrightarrow$ (h7) $\Longrightarrow$
(h8).}

 We keep the term $k+1$ in (h1)
 because it arises naturally in many examples.
 The first implication in Theorem 1
 was given in Stoyanov et al.\,(2014) recently, while the rest,
 more or less, are known in the
 literature.
 The Carleman quantity ${C}[F]$ in (h7) is calculated
 from all even order moments of $F$. Theorem 1 contains most checkable criteria for moment determinacy in the Hamburger
 case.

\noindent {\bf Remark 1.}
 Some other M-det
criteria exist in the literature, but they are seldom used. See, for example, (ha) and (hb) below:\\
{(h2)} $X$ has a mgf (Cram\'er's condition)\\
$\Longleftrightarrow$ (ha)
$\sum_{k=1}^{\infty}\frac{m_{2k}}{(2k)!}x^{2k}$ converges in an
interval $|x|<x_0$ (Chow and Teicher 1997,\ p.\,301)\\
$\Longrightarrow$ (hb) $\sum_{k=1}^{\infty}\frac{m_k}{k!}x^k$
converges
in an interval $|x|<x_0$ (Billingsley 1995, p.\,388)\\
$\Longrightarrow$ (h7)
${C}[F]=\sum_{k=1}^{\infty}m_{2k}^{-1/(2k)}=\infty$ (Carleman's
condition)\\
$\Longrightarrow$ (h8) $X$ is M-det on ${\mathbb R}.$\\It might look
strange that the convergence of subseries in the above (ha) implies
the convergence of the whole series in (hb), but remember that the
convergence in (ha) holds true for all $x$ in a neighborhood of
zero, not just for a fixed $x$. Billingsley (1995) proved the
implication that (hb) $\Longrightarrow$ (h8) by a version of
analytic continuation of characteristic function, but it is easy to
see that (hb) also implies (h7) and hence $X$ is M-det on ${\mathbb
R}.$

In Theorem 1, {Carleman's condition} (h7) is the weakest checkable
condition for
 $X$ to be M-det on ${\mathbb R}$. To prove Carleman's criterion that (h7) implies (h8), we may
 apply the approach of {quasi-analytic functions} (Carleman 1926, Koosis
 1988), or the approach of {L\'evy distance} (Klebanov and Mkrtchyan
 1980). For the latter, we recall the following result.\\
{\bf Klebanov and Mkrtchyan's (1980) Theorem.}  Let $F$ and $G$ be
two distribution functions on ${\mathbb R}$ and let their first $2n$
moments exist and coincide: $m_k(F) = m_k(G) = m_k$, $k =
1,2,\ldots, 2n\ (n\ge 2).$ Denote the sub-quantity
$C_n=\sum_{k=1}^nm_{2k}^{-1/(2k)}$. Then
$$L(F,G) \le c_2 \frac{\log(1+C_{n-1})}{(C_{n-1})^{1/4}},$$
where $L(F,G)$ is the L\'evy distance and $c_2 = c_2(m_2)$ depends
only on $m_2$. \\Therefore, {Carleman's condition} (h7) implies that
$F=G$ by letting $n\to\infty$ in Klebanov and Mkrtchyan's (1980)
Theorem.  It worths mentioning that Carleman's condition is
sufficient, but not necessary, for a distribution to be M-det. For
this, see Heyde (1963b), Stoyanov and Lin (2012, Remarks 5 and 7) or
Stoyanov (2013, Section 11).

On the other hand, the statement (h1) in Theorem 1 is the
strongest checkable condition for
 $X$ to be M-det on ${\mathbb R}$, which means that the growth rate of
 even order moments is less than or equal to two. The condition
 (h1) however has its advantage: for some cases, it is easy to estimate the
 growth rate (see the example below), because the common factors in the two even
 order moments, $m_{2(k+1)}$ and $m_{2k}$, can be cancelled out as $n$ tends to infinity.
\smallskip\\
\noindent{\bf Example 1.} Consider the double generalized gamma
random variable {$\xi \sim DGG(\alpha, \beta, \gamma)$} with density
function $f(x)=c|x|^{\gamma-1}\exp[{-\alpha |x|^{\beta}}],~x\in
{\mathbb R},$ where $\alpha, \beta, \gamma>0, f(0)=0$ if $\gamma\ne
1$, and $c=\beta\alpha^{\gamma/\beta}/(2\Gamma(\gamma/\beta))$ is
the norming constant. Then the $n$th power $\xi^n$ is M-det if $1\le
n\leq \beta$. To see this known result, we calculate the ratio of
even order moments of $\xi^n$:
\begin{eqnarray*}
\frac{{\bf
E}[\xi^{2n(k+1)}]}{{\bf E}[\xi^{2nk}]}=
\frac{\Gamma((\gamma+2n(k+1))/\beta)}{\alpha^{2n/\beta}\Gamma((\gamma+2nk)/\beta)}
\approx \left(2n/(\alpha\beta)\right)^{2n/\beta} {(k+1)^{2n/\beta}}~~
\hbox{as}~ k\rightarrow \infty,
\end{eqnarray*}
by using the approximation of the gamma function:
$\Gamma(x)\approx
\sqrt{2\pi}x^{x-1/2}e^{-x}~\hbox{as}~x\rightarrow \infty.$
 Therefore,  $\xi^n$ is M-det if $n\leq \beta$, by the criterion (h1).
 In fact, for odd integer $n\ge 1,$\
$\xi^n$ is M-det iff $n\le \beta,$ and for even integer $n\ge 2,$\
$\xi^n$ is M-det iff $n\le 2\beta,$ regardless of parameter
$\gamma.$ For further
 results about this distribution and its extensions, see Lin and Huang (1997),
 Pakes et al. (2001) and Pakes (2014, Theorem 8.3), as well as Examples 3 and 5 below.

\noindent {\bf Remark 2.}
 We give here a direct proof of
 the equivalence of statements  (h2), (h3) (h5) and (h6).
 First, for any nonnegative $X,$ we have the equivalence of the following four statements (to be shown later):\\
 ${\bf E}[e^{c\sqrt{{X}}}]<\infty\ \hbox{for some constant}\ ~ c>0$\\ {iff}
${m_k}\leq c_0^k(2k)!, ~k=1,2,\ldots,$ for some constant $c_0>0$\\
 {iff}
$\limsup_{k\rightarrow \infty}\frac{1}{k}\,{m_k}^{1/(2k)}<\infty$\\
{iff} ${m_{{k}}}={\cal O}(k^{2k})$ as $k\to\infty.$\\
Next, consider a general $X$ with ${\bf E}[e^{t{|X|}}]<\infty$ for
$0\le t<c$, namely, ${\bf E}[e^{t{\sqrt{{|X|^2}}}}]<\infty$ for some
constant $t>0$.  Then the $k$th moment of $|X|^2$ is exactly the
$2k$th moment of $X$ and we have immediately the following
equivalences
(by taking $|X|^2$ as the above nonnegative $X$): \\
(h2) $X$ has a mgf\\
{iff} (h6)
${m_{2k}}\leq c_0^k(2k)!, ~k=1,2,\ldots,$ for some constant $c_0>0$\\
 {iff}
$\limsup_{k\rightarrow
\infty}\frac{1}{k}\,{m_{2k}}^{1/(2k)}<\infty$\
(iff (h3) $\limsup_{k\rightarrow \infty}\frac{1}{2k}\,m_{2k}^{1/(2k)}<\infty$)\\
{iff} ${m_{{2k}}}={\cal O}(k^{2k})$ as $k\to\infty$\ (iff (h5)
$m_{{2k}}={\cal O}((2k)^{2k})$ as $k\to\infty$).

We now present the checkable M-det criteria in the
 Stieltjes case.  Consider $X\sim F$ on ${\mathbb R}_+=[0,\infty)$ with finite
$m_k=\mu_k={\bf E}[X^k]$ for all positive integers $k$, and define
the following
statements, in which `s' stands for `Stieltjes'.\\
(s1) $\frac{m_{k+1}}{m_{k}}={\cal O}((k+1)^2)={\cal O}(k^2)$ as $k\to\infty$.\\
(s2) $\sqrt{X}$ has a mgf (Hardy's condition), i.e., ${\bf
E}[e^{c\sqrt{X}}]<\infty$ for
some  constant $c>0$. \\
(s3) $\limsup_{k\to\infty}\frac{1}{k}m_{k}^{1/(2k)}<\infty.$ \\
(s4) $m_{{k}}={\cal O}(k^{2k})$ as $k\to\infty.$\\
(s5) $m_{{k}}\le c_0^k\,(2k)!,\ k=1,2,\ldots,$ for some constant
$c_0>0.$ \\
(s6) ${C}[F]=\sum_{k=1}^{\infty}m_{k}^{-1/(2k)}=\infty$ (Carleman's
condition).\\
(s7) $X$ is M-det on ${\mathbb R}_+.$

\noindent {\bf Theorem 2.} Under the above settings, if $X\sim F$ on
${\mathbb R}_+$ satisfies one of the conditions (s1) through (s6),
then $X$ is M-det on ${\mathbb R}_+$. Moreover, (s1) implies (s2),
(s2) through (s5) are equivalent, and (s5) implies (s6).
In other words, the following chain of implications holds:\\
\centerline {(s1) $\Longrightarrow$ (s2) $\Longleftrightarrow$ (s3)
$\Longleftrightarrow$ (s4) $\Longleftrightarrow$ (s5)
$\Longrightarrow$ (s6) $\Longrightarrow$ (s7).}

 The first implication above
 was given in Lin and Stoyanov (2015). Note that the moment conditions here are in terms of moments of all positive (integer) orders,
 rather than even order moments as in the Hamburger case. For example, the statement (s1) means that
the growth rate of all
  moments (not only for even order moments) is less than or equal to two. Like Theorem 1,
  Theorem 2 contains most  checkable criteria for moment determinacy in the Stieltjes
 case.
 Hardy (1917/1918) proved that
(s2) implies (s7) by two different approaches. Surprisingly, Hardy's
criterion has been ignored for about one century since publication.
The following new
 characteristic properties of (s2) are given in Stoyanov and Lin
 (2012),
 from which the equivalence
 of (s2) through (s5) follows immediately.
\smallskip\\
\noindent
{{\bf Lemma 1.}} Let $a$ be a positive constant and $X$ be a nonnegative random variable. \\
(i) If ${\bf E}[\exp({cX^a})]<\infty$ for some constant $c>0$,
then $m_k\leq \Gamma(k/a +1)c_0^k, ~k=1,2,\ldots,$ for
some constant $c_0>0.$ \\
(ii) Conversely, if, in addition, $a\leq 1$, and $m_k\leq
\Gamma(k/a +1)c_0^k, ~k=1,2,\ldots,$ for some constant $c_0>0,$
 then ${\bf E}[\exp({cX^a})]<\infty$ for some constant
$c>0$.\smallskip\\ {{\bf Corollary 1.}} Let $a\in (0, 1]$ and
$X\ge 0$. Then ${\bf E}[\exp({cX^a})]<\infty\ \hbox{for some
constant}\ ~ c>0$ {iff} $m_k\leq \Gamma(k/a +1)c_0^k,
~k=1,2,\ldots,$ for some constant $c_0>0.$
\smallskip\\
\noindent {\bf Lemma 2.}  Let $a$ be a positive constant and $X$ be
a nonnegative random variable. Then $\limsup_{k\rightarrow \infty}
\frac{1}{k}\,m_k^{a/k}$$<\infty$
 iff $m_k\leq \Gamma(k/a+1)\,c_0^k,~k=1,2,\ldots$, for some constant
 $c_0>0.$\medskip\\
 {{\bf Corollary 2.}} Let $a\in (0, 1]$ and $X\ge 0.$
 Then ${\bf E}[\exp({cX^a})]$$<\infty$ for some constant $c>0$\ \ iff\ \
$\limsup_{k\rightarrow \infty} \frac{1}{k}\,m_k^{a/k}$ $<\infty$.

We mention that for any nonnegative $X$, its mgf exists iff
$\limsup_{k\rightarrow \infty} \frac{1}{k}\,m_k^{1/k}$ $<\infty$ due
to Corollary 2. This in turn implies the equivalence of (h2) and
(h4) in Theorem 1 for the Hamburger case.  More general results in
terms of absolute moments are given below. For easy comparison, some
statements are repeated here.
\smallskip\\
\noindent {\bf Equivalence Theorem A} (Hamburger case){\bf .} Let
$p\ge 1$ be a constant and the random variable $X\sim F$ on
${\mathbb R}.$ Denote $m_k={\bf E}[X^k]$ for integer $k\ge 1$ and
let $\mu_{\ell}={\bf E}[|X|^{\ell}]<\infty$ for all $\ell>0$. Then
the following
statements are equivalent:\\
(a) $X$ satisfies Cram\'er's condition, namely, the moment
generating function of $X$ exists.\\
(b) $\mu_k\le c_0^kk!, k=1,2,\ldots,$ for some constant $c_0>0$.\\
(c) $\mu_{pk}\le c_0^k\Gamma(pk+1), k=1,2,\ldots,$ for some
constant $c_0>0$.\\
(d) $\limsup_{k\to \infty}\frac{1}{pk}\mu_{pk}^{1/(pk)}<\infty$.\\
(e) $m_{2k}\le c_0^k(2k)!, k=1,2,\ldots,$ for some
constant $c_0>0$.\\
(f) $\limsup_{k\to \infty}\frac{1}{2k}m_{2k}^{1/(2k)}<\infty$.\\
{\bf Proof.} The equivalence of (a), (b), (e) and (f) was given in
Theorem 1. To prove the remaining relations, denote $X_*=|X|$ and
write $Y_p=X_*^p$ and $\nu_{k,p}={\bf E}[Y_p^k]=\mu_{pk}$. Then note
further that ${\bf E}[e^{cX_*}]={\bf E}[e^{c(Y_p)^{1/p}}]<\infty$
for some constant $c>0$ iff $\nu_{k,p}\le c_0^k\Gamma(pk+1),
k=1,2,\ldots,$ for some constant $c_0>0$  (by taking $a=1/p$ and
$X=Y_p$ in Lemma 1) iff (c) holds true.  On the other hand,
$\nu_{k,p}\le c_0^k\Gamma(pk+1), k=1,2,\ldots,$ for some constant
$c_0>0$ iff $\limsup_{k\to
\infty}\frac{1}{k}\nu_{k,p}^{1/(pk)}<\infty$ (by Lemma 2) iff (d)
holds true. The proof is complete.

The above statements (e) and (f) are special cases of (c) and (d)
with $p=2$, respectively. Similarly, we give the following
equivalence theorem without proof for Stieltjes case.
\smallskip\\
\noindent {\bf Equivalence Theorem B} (Stieltjes case){\bf .} Let
$p\ge 1$ be a constant. Let the random variable $0\le X\sim F$ on
${\mathbb R}_+$ with finite $m_k=\mu_k={\bf E}[X^k]$ for all
integers $k\ge 1$. Then the following
statements are equivalent:\\
(a) $X$ satisfies Hardy's condition, namely, the moment
generating function of $\sqrt{X}$ exists.\\
(b) $\mu_k\le c_0^k(2k)!, k=1,2,\ldots,$ for some constant $c_0>0$.\\
(c) $\mu_{pk}\le c_0^k\Gamma(2pk+1), k=1,2,\ldots,$ for some
constant $c_0>0$.\\
(d) $\limsup_{k\to \infty}\frac{1}{pk}\mu_{pk}^{1/(2pk)}<\infty$.\\
(e) $\limsup_{k\to \infty}\frac{1}{k}\mu_{k}^{1/(2k)}<\infty$.
\medskip\\
\noindent{\bf 3. Checkable Criteria for Moment Indeterminacy}

 In
this section we consider the checkable criteria for moment
indeterminacy. In 1945, Krein proved the following remarkable
criterion in the {Hamburger case}.
\smallskip\\
\noindent{\bf Krein's Theorem.} Let $X\sim F$ on ${\mathbb R}$ have
a positive density function $f$ and finite moments of all positive
orders. Assume further that the Lebesgue logarithmic integral
\begin{eqnarray}{K}[f]\equiv\int_{-\infty}^{\infty}\frac{-\log f(x)}{1+x^2}dx<\infty.
\end{eqnarray}
Then $F$ is M-indet on ${\mathbb R}$.

We call the logarithmic integral $K[f]$ in (1) the Krein integral
for the density $f$. Graffi and Grecchi (1978)  as well as Slud
(1993) proved independently the counterpart of Krein's Theorem for
the {Stieltjes case} by the method of symmetrization of a
distribution on ${\mathbb R}_+$. To give a constructive and complete
proof, we however need Lemma 3 below (see, e.g., Lin 1997, Theorem
3, and Rao et al.\,2009, Remark 8).
\smallskip\\
\noindent{\bf  Graffi, Grecchi  and Slud's Theorem.} Let $X\sim F$
on ${\mathbb R}_+$ have a positive density function $f$ and finite
moments of all positive orders. Assume further that the integral
\begin{eqnarray}{K}[f]=\int_{0}^{\infty}\frac{-\log
f({x^2})}{1+x^2}dx<\infty. \end{eqnarray} Then $F$ is M-indet on
${\mathbb R}_+$ and hence M-indet on ${\mathbb R}.$
\smallskip\\
\noindent {\bf Lemma 3.}   Let $Y$ have a symmetric distribution $G$
with density $g$ and finite moments of all positive orders. If the
integral
$${K}[g]=\int_{-\infty}^{\infty}\frac{-\log g(x)}{1+x^2}dx<\infty,$$
then there exists a {\it symmetric} distribution $G_*\ne G$
having the same moment sequence as $G$.\\
{\bf Proof.} By the assumptions of the lemma,  there exists a
complex-valued function $\phi$ such that $|\phi|=g$ (in the sense of
almost everywhere) and
$$\int_{-\infty}^{\infty}\phi(x)e^{itx}dx=0,~~~t\ge 0$$
(see the proof of Theorem 1 in Lin 1997 for details, and Garnett
1981, p.\,66,
 for the construction
of $\phi$). The last equality implies that
$$\int_{-\infty}^{\infty}x^k\phi(x)e^{itx}dx=0,~~t\ge 0,~k=0,1,2,\ldots.$$
In particular,
$$\int_{-\infty}^{\infty}x^k\phi(x)dx=0,~~k=0,1,2,\ldots. $$
Let $\phi=\phi_1+i\phi_2$, then both $\phi_j$ are real and
$|\phi_j|\le g$. We have
$$\int_{-\infty}^{\infty}x^k\phi_j(x)dx=0,~j=1,2,~~k=0,1,2,\ldots.$$
\indent We split the rest of the proof into three cases: \\
\centerline{ (i)~ $\phi_1\ne 0,~\phi_2=0,$ ~~{(ii)}~ $\phi_1=
0,~\phi_2\ne 0,$~~{and}\ \ {(iii)}~ $\phi_1\ne 0,~\phi_2\ne 0.$} (i)
If $\phi_1$ is odd, then for each $t>0$, the function
$\phi_*(x):=\phi_1(x)\sin(tx)$ is even and
$$\int_{-\infty}^{\infty}x^k\phi_*(x)dx=0,~k=0,1,2,\ldots.$$
Take $g_*=g+\phi_*\ne g.$ Then $\phi\ge 0$ is even  and has the same
moment sequence as $g$. On the other hand, if $\phi_1$ is not odd,
then let first $\ell(x)=\frac{1}{2}[\phi_1(x)+\phi_1(-x)]$ which is
even and satisfies
$$\int_{-\infty}^{\infty}x^k\ell(x)dx=0,~k=0,1,2,\ldots.$$
Next, take $g_*=g+\ell\ne g,$ which  has the same moment sequence as $g$.\\
(ii) The proof of this case is similar to that of case (i).\\
(iii)  If one $\phi_j$ is not odd, then it is done as in (i) (by
taking $\ell(x)=\frac{1}{2}[\phi_j(x)+\phi_j(-x)]$ and
$g_*=g+\ell$). Suppose now that both $\phi_j$ are odd, then, by the
definition of $\phi,$  we further have
$\int_{-\infty}^{\infty}\phi(x)e^{itx}dx=0\ \forall~t\in {\mathbb
R}.$ Let $t>0$ be fixed and define the function
$\psi(x)=\phi_1(x)\sin(tx)+\phi_2(x)\cos(tx)$ (the imaginary part of
$\phi(x)e^{itx}$), then
$$\int_{-\infty}^{\infty}x^k\psi(x)dx=0,~\int_{-\infty}^{\infty}x^k
\psi(-x)dx=0,~k=0,1,2,\ldots.$$ Take
$m(x)=\frac{1}{2}[\psi(x)+\psi(-x)]=\phi_1(x)\sin(tx)\ne 0,$ which
is even and satisfies
$$\int_{-\infty}^{\infty}x^km(x)dx=0,~k=0,1,2,\ldots.$$
We have $g_*=g+m\ne g,$ which is nonnegative and has the same moment
sequence as $g$. The proof is complete.

 It should be noted that in the
logarithmic integral (2), the argument of the density function $f$
is $x^2$ rather than $x$ as in (1). Recently, Pedersen (1998)
improved Krein's Theorem by the concept of positive lower uniform
density sets and proved that it suffices to calculate the Krein
integral over the two-sided tail of the density function (instead of
the whole line).

 \noindent{\bf Theorem 3} (Pedersen 1998){\bf .} Let
$X\sim F$ on ${\mathbb R}$ have a density function $f$ and finite
moments of all positive orders. Assume further that the integral
\begin{eqnarray}{K}[f]=\int_{{|x|\ge c}}\frac{-\log f(x)}{1+x^2}dx<\infty \ \
\hbox{for some } c\ge 0.
\end{eqnarray}
Then $X$ is M-indet on ${\mathbb R}$.

See also H\"orfelt (2005) for Theorem 3 with a different proof
(provided by H.\,L.\, Pedersen).
 Pedersen (1998) also showed  by
giving an example that Krein's condition (1) is sufficient, but not
necessary, for a distribution to be M-indet. This corrected the
statement (2) in Leipnik (1981) about Krein's condition.
 On the other hand, Pakes (2001) and H\"orfelt (2005) pointed out the
counterpart of Pedersen's Theorem for the Stieltjes case. To prove
this result, we need Lemma 4 below.
\smallskip\\
\noindent{\bf Theorem 4} (Pakes 2001, H\"orfelt 2005){\bf .} Let
$X\sim F$ on ${\mathbb R}_+$ have a density function $f$ and finite
moments of all positive orders. Assume further that the integral
\begin{eqnarray}
{K}[f]=\int_{{x\ge c}}\frac{-\log f({x^2})}{1+x^2}dx<\infty \ \
\hbox{for some } c\ge 0.
\end{eqnarray}
Then $X$ is M-indet on ${\mathbb R}_+$ and hence M-indet on
${\mathbb R}.$
\smallskip\\
\noindent {\bf Lemma 4.} Let $0\le X\sim F$ with density $f$ and
finite moments of all positive orders. Let $Y\sim G$ with density
$g$ be the symmetrization of $\sqrt{X}$. If for some $c\ge 0$,
$${K}[g]=\int_{|x|\ge c}\frac{-\log g(x)}{1+x^2}dx<\infty,$$
then $X$ is M-indet on ${\mathbb R}_+$.\\
{\bf Proof.} Under the condition on the logarithmic integral of $g$,
Pedersen (1998, Theorem 2.2) proved that the set of polynomials is
not dense in L$^1({\mathbb R}, g(x)dx)$. This implies that
 the set of
polynomials is not dense in L$^2({\mathbb R}, g(x)dx)$ either (see,
e.g., Berg and Christensen 1981, or Goffman and Pedrick 2002,
p.\,162). Then proceeding along the same lines as in the proof of
Corollary 1 in Slud (1993), we conclude that the set of polynomials
is not dense in L$^2({\mathbb R}, f(x)dx)$. Therefore, $X$ is
M-indet on ${\mathbb R}$, which in turn implies that $X$ is
 M-indet on ${\mathbb R}_+$ due to Chihara's (1968) result in Fact A above. The proof is complete.

  Conversely, once we prove Theorem 4,  we can extend Lemma 4 as follows.
\smallskip\\
\noindent {\bf Lemma 4$^*$.} If $X\sim F$ on ${\mathbb R}$
 satisfies the conditions in Theorem 3, then $X^2$  is M-indet.\\
 {\bf Proof.} Apply
Theorem 4 above and Pakes et al.'s (2001) Theorem 3(i): If $X\sim F$
on ${\mathbb R}$ satisfies condition (3), then the Krein integral
$K[f_2]$ in (4) of $X^2$ is finite, where $f_2$ is the density of
$X^2$.

For the M-det case, a trivial analogue of Lemma 4$^*$  is the
following.

\noindent {\bf Lemma 4$^{**}$.} If $X\sim F$ on ${\mathbb R}$
 satisfies Carleman's condition (h7), then $X^2$ satisfies Carleman's condition (s6) and is M-det on ${\mathbb R}_+.$

For simplicity, all the conditions (1) through (4) are called
{Krein's condition}. For illustration of how to use Krein's and
Hardy's criteria, we now recover Berg's (1988)  results using these
powerful criteria (see also Prohorov and Rozanov 1969, p.\,167,
Pakes and Khattree 1992, Lin and Huang 1997, and Stoyanov 2000).
\smallskip\\
\noindent
{\bf Example 2.} Let $X$ have a {normal distribution} and $\alpha>0.$ Then\\
{(i) the odd power $X^{2n+1}$ is M-indet if $n\ge 1,$} and\\
{(ii) $|X|^{\alpha}$ is M-det iff $\alpha\le 4.$}\\
Without loss of generality, we assume that $X$ has a density
$f(x)=\frac{1}{\sqrt{\pi}}\exp({-x^2}),\ x\in {\mathbb R},$ namely,
$\sqrt{2}X$\, has a standard normal distribution. We discuss these results in three steps.\\
(I) Berg (1988) proved the moment indeterminacy of distributions by
giving examples.  For part (i), he calculated first the density of
$X^{2n+1}$:
$$  f_n(x)=\frac{1}{(2n+1)\sqrt{\pi}}|x|^{-2n/(2n+1)}\exp({-{|x|^{2/(2n+1)}}}),\
x\in {\mathbb R},$$ and then constructed the density function
$$f_{r,n}(x)=f_n(x)\{1+r{[\cos(\beta_n|x|^{2/(2n+1)})-
\gamma_n\sin(\beta_n|x|^{2/(2n+1)})]}\}\equiv f_n(x)\{1+rp_n(x)\},\ x\in {\mathbb R},$$
where $|r|\le \sin\frac{\pi}{2(2n+1)},\
\beta_n=\tan\frac{\pi}{2n+1}$ and $\gamma_{\alpha}=
\cot\frac{\pi}{2(2n+1)}.$  It is seen that $f_{r,n}\ne f_n$ if $r\ne
0$ and $n\ge 1,$ but $f_{r,n}$ and $f_n$ have the same moment
sequence
 because the product of the density $f_n$ and the function $p_n$ defined above has
 vanishing moments by a tedious calculation:
$$\int_{-\infty}^{\infty}x^kf_n(x){p_n(x)}dx=0,\ k=0,1,2,\ldots.$$ This
proves part (i). Alternatively, we note however that the Krein
integral
$${{K}[f_n]=\int_{-\infty}^{\infty}\frac{-\log
f_n(x)}{1+x^2}dx=C+\int_{-\infty}^{\infty}\frac{|x|^{2/(2n+1)}}{1+x^2}dx<\infty\
(\hbox{if}\ n\ge 1),}$$ which implies by Krein's Theorem that the
odd power
$X^{2n+1}$ is M-indet if $n\ge 1.$\\
(II)  For part (ii), the density of $|X|^{\alpha}$ is
$$ f_{\alpha}(x)=\frac{2}{\alpha\sqrt{\pi}}x^{1/\alpha-1}
\exp({-{x^{2/\alpha}}}),\ x\ge 0.$$ If $\alpha>4,$ Berg
constructed again the density function
$$f_{r,\alpha}(x)=f_{\alpha}(x)\{1+r{[\cos(\beta_{\alpha}x^{2/{\alpha}})-
\gamma_{\alpha}\sin(\beta_{\alpha}x^{2/{\alpha}})]}\}\equiv
f_{\alpha}(x)\{1+rp_{\alpha}(x)\},\ x\ge 0,$$ where $|r|\le
\sin(\pi/{\alpha}),\ \beta_{\alpha}=\tan(2\pi/{\alpha})$ and
$\gamma_{\alpha}= \cot(\pi/{\alpha}).$  Then $f_{r,{\alpha}}\ne
f_{\alpha}$ if $r\ne 0$ and ${\alpha}> 4,$ but $f_{r,{\alpha}}$ and
$f_{\alpha}$ have the same moment sequence
 because
$$\int_{0}^{\infty}x^kf_{\alpha}(x){p_{\alpha}(x)}dx=0,\
k=0,1,2,\ldots.$$ Therefore, $|X|^{\alpha}$ is M-indet if $\alpha>
4.$ Again, we see that the Krein integral (in Stieltjes case)
$${{K}[f_{\alpha}]=\int_{0}^{\infty}\frac{-\log
f_{\alpha}(x^2)}{1+x^2}dx=C+\int_{0}^{\infty}\frac{x^{4/{\alpha}}}{1+x^2}dx<\infty\
(\hbox{if}\ {\alpha}> 4).}$$ So the required result follows
immediately
from Krein's criterion (4). \\
(III) For the rest of part (ii), Berg calculated the $k$th moment
of $|X|^{\alpha}$:
$$m_{\alpha,k}=\int_0^{\infty}x^kf_{\alpha}(x)dx=\frac{1}{\sqrt{\pi}}
\Gamma\left(\frac{\alpha k+1}{2}\right),\ k=0,1,2,\ldots.$$ By
Stirling's formula, $m_{\alpha,k}^{1/k}\approx ck^{\alpha/2}$ as
$k\to \infty,$ and hence the Carleman quantity (in Stieltjes case)
is equal to
$${C}[f_{\alpha}]=\sum_{k=1}^{\infty}m_{\alpha,k}^{-1/(2k)}=\infty\
(\hbox{if}\ {\alpha}\le 4).$$ This proves the necessary part of
(ii). Instead, we note that if $\alpha\in(0, 2],$ the mgf of
$|X|^{\alpha}$ exists by its density function above, and hence
$|X|^{2\alpha}$ is M-det by Hardy's criterion.

There are some ramifications of the moment problem for normal
random variables. For example, Slud (1993) investigated the moment
problem for polynomial forms in normal random variables, while
 H\"orfelt (2005) studied the moment
problem for some Wiener functionals which extend Berg's results in
Example 2. Besides, Lin and Huang (1997) treated the double
generalized Gamma (DGG) distribution as an extension of the normal
one and found the necessary and sufficient conditions for powers of
DGG random variable to be M-det.
\medskip\\
\noindent{\bf 4. Stieltjes Classes for M-indet Distributions}

 {Stieltjes (1894)} observed that some positive measures, e.g., $\mu(dx)=
e^{-x^{1/4}}dx$ or  $x^{n-\log x}dx$ ($n$ is an integer), are not
unique by moments. This might be the starting point of
T.\,J.\,Stieltjes to study the moment problem (see Kjeldsen 1993).
It was C.\,C.\,Heyde who first presented this phenomenon in
probability language and proved in 1963 that the {lognormal
distribution} is M-indet by giving the example described next.
Consider the standard lognormal density
$$f(x)=\frac{1}{\sqrt{2\pi}}x^{-1}\exp[-\frac{1}{2}(\log x)^2],\
x>0,$$ with moment sequence $\{\exp(k^2/2)\}_{k=1}^{\infty}.$ Then, for each
$\varepsilon\in[-1,1]$,
\begin{eqnarray*}\int_0^{\infty}{f(x)}[1+\varepsilon\sin(2\pi\log x)]x^kdx
=\int_0^{\infty}{f(x)}x^kdx \ \ \forall\
k=0,1,2,\ldots\end{eqnarray*} because the product of the density
$f(x)$ and the function $\sin (2\pi \log x)$ has vanishing moments:
$$\int_0^{\infty}
f(x)[\sin (2\pi \log
x)]x^kdx=\frac{e^{k^2/2}}{\sqrt{2\pi}}\int_{-\infty}^{\infty}e^{-x^2/2}\sin(2\pi
x)dx=0 \ \ \forall\ k=0,1,2,\ldots.$$ \indent There are many other
distributions  having the same moment sequence as the above
lognormal with mean $\sqrt{e}$, including (i) the ones with density
$f$ satisfying the functional equation: $f(qx)=q^{-1/2}xf(x),$ where
$q=1/e\in (0,1)$ (see, e.g.,  L\'opez-Garc\'ia 2011, Theorem 1), or,
more generally, (ii) the distributions $F$ satisfying
$F(x)=e^{-1/2}\int_0^{ex}udF(u),\ x\ge 0$ (Pakes 1996, Section 3).
The latter showed that each such $F$ corresponds to a finite measure
in the interval $(1/e,1]$ and vice versa. Hence the cardinality of
the set of all solutions to the functional equation is
${\aleph}_2=2^{\mathbb R}.$ All these distributions are called the
solutions to the lognormal moment problem
 (see also Chihara 1970,  Leipnik 1982, Pakes  2007 and
Christiansen 2003).

Recently, Stoyanov (2004) formulated Stieltjes classes for M-indet
absolutely continuous distributions as follows.  Let $X\sim F$ have
an M-indet distribution on ${\mathbb R}$ with density $f$. A
{Stieltjes class ${\cal S}$ for $F$} is defined by
$${\cal S}={\cal S}(f, p)=\{f_{\varepsilon}:
f_{\varepsilon}(x)=f(x)[1+\varepsilon p(x)],\ x\in {\mathbb R},\
\varepsilon\in [-1,1]\},$$ where $p$ is a measurable function
(called a perturbation function) such that $|p(x)|\le 1$ and
$${\int_{-\infty}^{\infty}f(x)p(x)x^kdx=0},\ \ k=0,1,2,\ldots.$$

We note that for a given M-indet distribution, the choice of the
function $p$ might not be unique. Besides the previous lognormal and
normal results of Heyde (1963) and Berg (1988),  some other
perturbation functions are given below:
\begin{enumerate}

\item If $X$ has a generalized Weibull density
$f(x)=\frac{1}{24}\exp(-x^{1/4}),\,x>0$, then
$p(x)=\sin(x^{1/4}),\, x>0$ (Stieltjes 1894, Serfling 1980).

\item If $X$ has a density function $f(x)=c x^{-\log x},$ $x>0$,
where $c$ is a norming constant, then
 we choose
$p(x)=\sin(2\pi\log x),\,x>0$ (Stieltjes 1894).

\item If $X$ has a gamma density with parameter $\alpha>0$, then $X^{\beta}$ is M-indet provided
$\beta>\max\{2,2\alpha\}$, and we can choose
$$p(x)=\sin(\alpha\pi/\beta)[\cos(\tan(\pi/\beta)x^{1/\beta})-
\cot(\alpha\pi/\beta)\sin(\tan(\pi/\beta)x^{1/\beta})],\ x>0$$
(Targhetta 1990).

\item  If $X$ has a density function $f(x)=c \exp(-\alpha |x|^{\rho}),$ $x\in {\mathbb R}$,
where $\alpha>0,\, \rho\in(0,1)$ and $c$ is a norming constant,
then we choose $p(x)=\cos(\alpha |x|^{\rho}),\,x\in {\mathbb R}$
(Prohorov
 and Rozanov 1969, p.\,167).

\item For the log-skew-normal distribution with parameter $\lambda>0$,
we choose the  perturbation function
 $$ p(x)=\frac{\ell(x-1)}{\ell(x)}\frac{\sin[\pi\log(x-1)]}{\Phi(\lambda\log
x)},\ \ \hbox{if}\ x>1,$$ and $p(x)=0,$ otherwise,
 where $\ell$ is the density of
standard lognormal LN(0,1) and $\Phi$ is the standard normal
distribution (Lin and Stoyanov 2009).
\end{enumerate}

Several systematic approaches for constructing Stieltjes classes are
available. For example, for any M-indet
distribution $F$ on $(0,\infty)$ with density $f$ bounded from below as\\
\centerline {$f(x)\ge A\exp(-\alpha x^{\beta}),$ $x>0$, where $A>0,\
\alpha>0$ and $\beta\in (0, 1/2)$ are constants,} we find first a
complex-valued function $g$
satisfying\\
(i) $g$ is analytic in ${\mathbb C}_+\setminus \{0\},$ where
${\mathbb C}_+=\{z:
\hbox{Im}\,z\ge 0\}$ is the upper half-plane,  and\\
(ii) $g(x)\in{\mathbb R}$, $x>0,$ and\ \ $|g(z)|\le
A\exp(-\alpha |z|^{\beta})$, $z\in {\mathbb C}_+\setminus \{0\}.$\\
Then choose the perturbation function
$p(x)=[\hbox{Im}\,g(-x)]/f(x),\ x>0$ (Ostrovska 2014).

 On the other
hand, given any Stieltjes class ${\cal S}(f, p)$ defined above and a
positive random variable $V$ with distribution $H$ and finite
moments of all positive orders, we can construct a new Stieltjes
class ${\cal S}(f^*, p^*)$ by random scaling:
$Y_{\varepsilon}:=VX_{\varepsilon},$ $\varepsilon\in [-1,1],$ where
 the random variable $X_{\varepsilon}$ has density
$f_{\varepsilon},$  $V$ is independent of $X_{\varepsilon},$
$f^*=f^*_0$ is the density of $Y_0=VX,$ and  the perturbation
function $p^*$ satisfies
$f^*(x)p^*(x)=\int_0^{\infty}v^{-1}f({x}/{v})p({x}/{v})dH(v),\ x\in
{\mathbb R}$ (Pakes 2007, Section 5).

For more  perturbation functions, see Stoyanov (2004), Stoyanov and
Tolmatz (2004, 2005), Ostrovska and Stoyanov (2005), G\'omez and
L\'opez-Garc\'ia (2007), Penson et al.\,(2010), Wang (2012), Kleiber
(2013, 2014) and Ostrovska (2016).\bigskip\\
{\bf 5. Converse Criteria}

 In this section we present some converses
to the previous M-(in)det criteria.  Recall that for the Stieltjes
case, if (s1) above holds true, i.e., ${m_{k+1}}/{m_{k}}={\cal
O}((k+1)^2)$ as $k\to\infty$, then $X$ is M-det on ${\mathbb R}_+.$
One might guess that if the moments $\{m_{k}\}_{k=1}^{\infty}$ grow
faster, then $X$ becomes M-indet. This is true under one more
condition defined below (see Lin 1997, Stirzaker 2015, p.\,223,
Kopanov and Stoyanov 2017, or Stoyanov and Kopanov 2017).

\noindent{\bf Condition L:} Suppose, in the Stieltjes case,  that
$f$ is a density function on ${\mathbb R}_+$ such that for some
fixed $x_0\ge 0,$\ $f$ is strictly positive and differentiable on
$(x_0,\infty)$ and
$$
L_{f}(x):={-\frac{xf^{\prime}(x)}{f(x)}\nearrow
\infty~~\hbox{as}~~x_0<x\rightarrow \infty.}
$$
In the Hamburger case we require the density $f(x), \ x \in {\mathbb
R},$
to be symmetric about zero.\smallskip\\
\noindent{\bf Theorem 5.} Let $X$ be a nonnegative random variable
with distribution $F$ and let its moments grow fast in the sense
that $m_{k+1}/m_k \geq c(k+1)^{{2+\varepsilon}}$ for all large $k$,
where $c$ and $\varepsilon$ are positive constants.
 Assume further that $X$
has a density function $f$ which satisfies {Condition L}. Then $X$
is {M-indet}.

Note that in the above theorem, $X$ is M-indet on ${\mathbb R}_+$
iff it is M-indet on ${\mathbb R}$ because  $X$ has a density.
For the Hamburger case, we have the following.\smallskip\\
\noindent{\bf Theorem 6.} Suppose the moments of $X\sim F$ on
${\mathbb R}$ grow fast in the sense that $m_{2(k+1)}/m_{2k}\ge
 c(k+1)^{{2+\varepsilon}}$ for all large $k$, where $c$ and
$\varepsilon$ are positive constants. Assume further that $X$ has a
density function $f$ which is symmetric about zero and satisfies
{Condition L}. Then $X$ satisfies Krein's condition, and hence both
$X$ and $X^2$ are {M-indet}.

The crucial point in the proofs of Theorems 5 and 6 is to prove that
the Krein integral ${K}[f]<\infty$ by Condition L and the moment
condition. The M-indet property of $X^2$ in Theorem 6 is due to
Lemma 4$^*$ and Fact A above. Similarly, we have the following
results for other criteria (s4) and (h5) (see Lin and
Stoyanov 2015, 2016, and Stoyanov et al.\,2014).\smallskip\\
\noindent {\bf Theorem 7} (Stieltjes case){\bf .} Let $X\sim F$ on
${\mathbb R}_+$ and let its moments grow fast in the sense that
 $m_k \ge
c\,k^{({2+\varepsilon})k}, \ k=1,2,\ldots,$ for some positive
constants $c$ and $\varepsilon$. Assume further that $X$ has a
density function $f$ which satisfies {Condition L}. Then $X$ is
{M-indet}.\smallskip\\ \noindent {\bf Theorem 8} (Hamburger
case){\bf .} Suppose the moments of $X\sim F$ grow fast in the sense
that $m_{2k} \ge c(2k)^{{(2+\varepsilon)}k},$  $k=1,2,\ldots,$ for
some positive constants $c$ and $\varepsilon.$ Assume further that
$X$ has a density function $f$ which is symmetric about zero and
satisfies {Condition L}. Then $X$ satisfies Krein's condition, and
hence both $X$ and $X^2$ are {M-indet}.

Note that Condition L also applies to converse M-indet criteria.
Actually, this is the original purpose of the condition, under which
${K}[f]=\infty$ implies ${C}[F]=\infty$ (Lin 1997).
The M-det property of $X^2$ in the next result is due to Lemma 4$^{**}$ and Fact A above.\smallskip\\
\noindent{\bf Theorem 9.} In Theorem 3 (Hamburger case), if the
Krein integral {${K}[f]=\infty$} and if $f$ satisfies {Condition L},
then $X$ satisfies Carleman's condition, and hence both $X$  and
$X^2$ are {M-det}.
\smallskip\\
\noindent{\bf Theorem 10.} In Theorem 4 (Stieltjes case), if the
Krein integral {${K}[f]=\infty$} and if $f$ satisfies {Condition L},
then $X$ satisfies Carleman's condition and is {M-det}.

\noindent {\bf Remark 3.} In view of Theorems 9 and 10 above, we
know that in the class of absolutely continuous distributions with
density functions satisfying Condition L, Krein's condition ((3) or
(4)) becomes necessary and sufficient for a distribution to be
M-indet.

\noindent {\bf Remark 4.} In the above converse results, it is
possible to replace Condition L by other  slightly weaker conditions
(mathematically) like those in Pakes (2001) and Gut (2002), but as
mentioned before, we focus only on the {\it checkable conditions} in
this survey. Interestingly, Condition L is closely related to a
useful concept in reliability theory. More precisely, if a
nonnegative random variable $X$ with density $F^{\prime}=f$
satisfies Condition L on ${\mathbb R}_+$ with $x_0=0$, then it has
an increasing generalized failure rate (by Theorem 1 in Lariviere
2006), namely, the product function $xf(x)/\overline{F}(x)$ (of $x$
and the failure rate) increases in $x.$

In addition to the previous problems for normal distributions, we
mention here some more variants for general cases, but we are not
going to pursuit all the moment problems.  To solve these problems,
we need to derive new auxiliary tools  case by case (like Lemma 5
below). Lin and Stoyanov (2002) and Gut (2003) studied the moment
problem for random sums of independently identically distributed
(i.i.d.) random variables. Stoyanov et al.\,(2014) and Lin and
Stoyanov (2015) investigated the moment problem for products of
i.i.d.\,random variables. In the next section we review the recent
results about products of independent random variables with {\it
different} distributions; for details, see Lin and Stoyanov (2016).
\medskip\\
\noindent{\bf 6. Moment Problem for Products of Random Variables}

Products of random variables occur naturally in stochastic modelling
of complex random phenomena in areas such as statistical physics,
quantum theory, communication theory, reliability theory and
financial modelling; especially in modern communications (see, e.g.,
Chen et al.\,2012, Springer 1979, and Galambos and Simonelli 2004).
We split the problem in question into three cases: (a) products of
nonnegative random variables, (b) products of random variables
taking values in ${\mathbb R}$, and (c) the mixed case. Moreover,
all random variables considered have finite moments of all positive
orders.
\medskip\\
{\bf 6.1. Products of Nonnegative Random Variables}

The M-det result (Theorem 11 below) is an easy consequence of
Theorem 2, while the hard part is the M-indet result (Theorem 12)
whose proof needs a delicate analysis.
\smallskip\\
\noindent{\bf Theorem 11.} {Let  $\xi_1, \ldots, \xi_n$ be
independent nonnegative random variables and let the moments
$m_{i,k}={\bf E}[\xi_i^k], \ i=1,\ldots,n,$
 satisfy the conditions:
\[
m_{i,k} = {\cal O}(k^{{a_i}k}) \ \mbox{ as } \ k \to \infty, \
\mbox{ for } \ i=1, \ldots, n,
\]
where $a_1, \ldots, a_n$ are positive constants. If the parameters
$a_i$ are such that $\sum_{i=1}^na_i\le 2,$ then the product
$Z_n=\Pi_{i=1}^n\xi_i$ satisfies Hardy's condition and is {M-det}.}
\smallskip\\
\noindent{\bf Theorem 12.} {Consider $n$ independent nonnegative
random variables, {$\xi_i \sim F_i,~i=1,2,\ldots, n$}, where $n \ge
2.$ Suppose that each {$F_i$ is absolutely continuous and has a
positive density $f_i$} on $(0,\infty)$ and
that the following conditions are satisfied: \\
(i) At least {one} of the densities $f_1(x),\ldots,f_n(x)$
is decreasing in $[x_0,\infty),$ where $x_0\ge 1$ is a constant.\\
(ii) For each $i=1,2,\dots, n,$ there exists  a constant $A_i>0$
such that the density $f_i$ and the tail function $\overline
{F_i}(x)=1-F_i(x)=\Pr(\xi_i>x)$ together satisfy the relation
\begin{eqnarray}
{f_i(x)/\overline{F_i}(x)\geq A_i/x~~\hbox{for}~~x\geq x_0,}
\end{eqnarray}
and there exist constants $B_i>0,~ \alpha_i>0,$ $\beta_i>0$ and
real $\gamma_i$ such that
\begin{eqnarray}
{\overline{F_i}(x)\geq B_ix^{\gamma_i}\exp({-\alpha_i
x^{\beta_i}})~~ \hbox{for}~~ x\geq x_0.}
\end{eqnarray}
If, in addition to the above, {$\sum_{i=1}^n1/{\beta_i}>2,$} then
the product $Z_n=\Pi_{i=1}^n\xi_i$ is {M-indet}.}

Let us explain the above conditions. In terms of reliability
language, the failure rate in (5) and the survival function in (6)
cannot approach zero too quickly. In other words, (5) and (6)
control the tail (decreasing) behavior of the related distributions
in some sense. There are three key steps in the proof of Theorem 12:
(i) represent the density function of the product $Z_n$ in multiple
integral form, (ii) estimate the lower bound of the density function
by truncating the two tails of this integral, and (iii) apply
Krein's criterion for the Stieltjes case. For estimation in the step
(ii), we need the following auxiliary tool which can be proved using
 integration by parts.
\smallskip\\
\noindent {\bf Lemma 5.} {Let $F$ be a distribution on ${\mathbb R}$
such that  (i) it has density $f$ on the subset $[a, ra],$ where
$a>0$ and $r>1,$ and  (ii) for some constant $A>0,$
${f(x)}/{\overline{F}(x)}\ge {A}/{x}$ on $[a,ra].$ Then}
$$\int_a^{ra}\frac{f(x)}{x}dx\ge
\left(1-\frac{1}{r}\right)\frac{A}{1+A}\frac{\overline{F}(a)}{a}.$$
\noindent {\bf Example 3.}  For illustration of how to use Theorems
11 and 12, consider the  {generalized gamma distributions}. We say
that $\xi \sim GG(\alpha, \beta, \gamma)$ if its density is of the
form
\begin{eqnarray*}
{f(x)=c\,x^{\gamma-1}\exp({-\alpha x^{\beta}}),\ ~x\geq 0.}
\end{eqnarray*}
Here $\alpha, \beta, \gamma>0$, $f(0)=0$ if $\gamma\ne 1,$ and
$c=\beta\alpha^{\gamma/\beta}/\Gamma(\gamma/\beta)$ is the norming
constant. Then we have the following characterization result (see
also Pakes 2014 for a much more general result with different proof):\\
Suppose that $\xi_1, \ldots, \xi_n$ are $n$ independent random
variables and let $\xi_i\sim {GG(\alpha_i,\beta_i,\gamma_i)},$
$i=1,\ldots,n.$
 Then the product $Z_n=\Pi_{i=1}^n\xi_i$ is M-det iff\,
 $\sum_{i=1}^n{1}/{\beta_i}$ $\leq 2.$
\smallskip\\
{\bf Example 4.} Consider the class of {inverse Gaussian
distributions}. We say that $\xi\sim IG(\mu,\lambda)$ if its density
is of the form
\begin{eqnarray*}
{f(x)=\left(\frac{\lambda}{2\pi
x^3}\right)^{1/2}\exp\left[-\frac{\lambda(x-\mu)^2}{2\mu^2x}\right],\
~x>0,}
\end{eqnarray*}
where $\mu, \lambda>0$ and $f(0)=0.$ It can be shown that the
product of two independent random variables is M-det if each one is
exponential or inverse Gaussian, while the product of three such
random variables is M-indet. For the powers of such random variables
and others, see, e.g., Lin and Huang (1997), Stoyanov (1999), Pakes
et al.\,(2001), Stoyanov et al.\,(2014) and Lin and Stoyanov (2015).
Here are some recent results.\\
 Let $\xi_1\sim IG(\mu_1,\lambda_1),\ \xi_2\sim
IG(\mu_2,\lambda_2)$ and $\eta\sim Exp(1)=GG(1,1,1)$ be three
independent random variables. Then both the products $\xi_1\eta$ and
$\xi_1\xi_2$ are M-det, while $\xi_1\xi_2\eta$ is M-indet.
\newpage
\noindent{\bf 6.2. Products of Random Variables Taking Values in}
${\mathbb
R}$\\
\indent For this Hamburger case, we have the counterparts of
Theorems 11 and 12 as follows.  In the proof of Theorem 14, the
symmetric condition on the densities plays a crucial role.
\smallskip\\
\noindent {\bf Theorem 13.}  {Let  $\xi_1, \ldots, \xi_n$ be
independent random variables and let the even order moments
$m_{i,2k}={\bf E}[\xi_i^{2k}], \ i=1,\ldots,n,$
 satisfy the conditions:
\[
m_{i,2k} = {\cal O}((2k)^{{2a_i}k}) \ \mbox{ as } \ k \to \infty, \
\mbox{ for } \ i=1, \ldots, n,
\]
where $a_1, \ldots, a_n$ are positive constants. If the parameters
$a_i$ are such that $\sum_{i=1}^na_i\le 1,$ then the product
$Z_n=\Pi_{i=1}^n\xi_i$ satisfies Cram\'er's condition and is
{M-det}. }\smallskip\\
{\bf Theorem 14.} {Consider $n$ independent random variables
{$\xi_i \sim F_i,$ $i=1,\ldots,$ $n,$}
 where $n \ge 2.$  Suppose each {$F_i$ has a positive density
 $f_i$ on ${\mathbb R}$ and
{symmetric about $0$}.} Assume further
that\\
(i) at least one of the densities
$f_1(x),\ldots,f_n(x)$ is decreasing in $[x_0,\infty),$ where $x_0\ge 1$ is a constant, and \\
(ii) for all $i$, $f_i/\overline{F_i}$ satisfies the condition (5):
$f_i(x)/\overline{F_i}(x)\geq A_i/x~~\hbox{for}~~x\geq x_0$, and
$\overline{F_i}$ satisfies the condition (6): $\overline{F_i}(x)\geq
B_ix^{\gamma_i}\exp({-\alpha_i
x^{\beta_i}})~~ \hbox{for}~~ x\geq x_0$.\\
If, in addition to the above, {$\sum_{i=1}^n1/{\beta_i}>1,$} then
the product $Z_n=\Pi_{i=1}^n\xi_i$ satisfies Krein's condition, and
hence both $Z_n$ and $Z_n^2$ are {M-indet}.}
\smallskip\\
{\bf Example 5.} Applying Theorems 13 and 14 to the product of
double generalized gamma random variables $\xi \sim DGG(\alpha,
\beta, \gamma),$ defined above, yields the following interesting result:\\
 {Suppose that $\xi_1, \ldots, \xi_n$ are $n$
independent random variables, and let $\xi_i\sim
{DGG(\alpha_i,\beta_i,\gamma_i)}, \ i=1,2,\ldots,n.$ Then the
product $Z_n=\Pi_{i=1}^n\xi_i$ is {M-det iff\,
 $\sum_{i=1}^n{1}/{\beta_i}\leq 1$}} iff\
 $Z_n^2$ is M-det.
\medskip\\
{\bf 6.3. The Mixed Case}

Finally, we consider the products of both types of random variables,
nonnegative and real ones taking values in ${\mathbb R}$. Recall
that this is the Hamburger case and the M-det criterion is similar
to Theorem 13 and omitted. The next result about an M-indet
criterion extends slightly Theorem 5.1 of Lin and Stoyanov (2016).
The proof is similar and is therefore omitted.
\smallskip\\ {\bf Theorem 15.} {Consider $n$ independent random
variables divided into two groups.
 The first group, {$\xi_1, \ldots,
 \xi_{n_0},$}
 consists of nonnegative variables, while all the variables
in the second group,  {$\xi_{n_0+1}, \ldots, \xi_n,$} take values in
${\mathbb R},$ where $1\le n_0<n.$ Suppose that {each $\xi_i \sim
F_i$ has a density $f_i$} and that $f_i,$\, $i=1,\ldots,n_0,$ are
positive on \,$(0,\infty),$ while {$f_j,$\, $j=n_0+1, \ldots, n,$
are  positive on ${\mathbb R}$ and
symmetric about 0.} Moreover, assume further that\\
(i) at least one  of the densities $f_j(x),$ $j=1,2,\ldots,n,$
is decreasing in $[x_0,\infty),$ where $x_0\ge 1$ is a constant, and \\
(ii) for all $i$, $f_i/\overline{F_i}$ satisfies the condition (5):
$f_i(x)/\overline{F_i}(x)\geq A_i/x~~\hbox{for}~~x\geq x_0$, and
$\overline{F_i}$ satisfies the condition (6): $\overline{F_i}(x)\geq
B_ix^{\gamma_i}\exp({-\alpha_i x^{\beta_i}})~~ \hbox{for}~~ x\geq
x_0$.\\
If, in addition to the above, {$\sum_{i=1}^n1/{\beta_i}>1,$} then
the product $Z_n=\Pi_{i=1}^n\xi_i$ satisfies Krein's condition, and
hence both $Z_n$ and $Z_n^2$ are {M-indet}.}

An application of the above theorem leads to the following
interesting result:\\
 The product of two independent random
variables  and its square are both M-indet if one random variable is
normal and the other is exponential, or chi-square, or inverse
Gaussian.
\medskip\\
\noindent {\bf Acknowledgments.} The author would like to thank the
Editor and two Referees for helpful comments and suggestions.
Especially, one Referee pointed out the result in Lemma 4$^*$. The
paper was presented at (1) the International Waseda Symposium,
February 29 -- March 3, 2016, held by Waseda University (Japan) and
(2) the second International Conference on Statistical Distributions
and Applications (ICOSDA), October 14--16, 2016, Niagara Falls, held
by Central Michigan University (USA) and Brock University (Canada).
The author thanks the organizers (1) Professor Masanobu Taniguchi
and (2) Professors Felix Famoye, Carl Lee and Ejaz Ahmed for their
kind invitations. The comments and suggestions of Professor Murad
Taqqu and other audiences are also appreciated.
\bigskip\\
\noindent {\bf References}
\begin{description}

\item {Akhiezer, NI}: {The Classical Problem of Moments
and Some Related Questions of Analysis}. Oliver $\&$ Boyd,
Edinburgh (1965) [Original Russian edition: Nauka,  Moscow
(1961)]

\item {Berg, C}: The cube of a normal distribution is
indeterminate. {Ann. Probab.} {\bf 16}, 910--913 (1988)

\item {Berg, C}: From discrete to absolutely continuous
solutions of indeterminate moment problems. {Arab. J. Math.
Sci.} {\bf 4}, 1--18 (1998)

\item {Berg, C, Chen, Y, Ismail, MEH}: Small eigenvalues of large Hankel matrices: the indeterminate case.
Math. Scand. {\bf 91}, 67--81 (2002)

\item Berg, C, Christensen, JPR: Density questions in
  the classical theory of moments. {Ann. Inst. Fourier
  (Grenoble)} {\bf 31}, 99--114 (1981)

\item {Billingsley, P}: {Probability and Measures}. 3rd edn.
 Wiley, New York (1995)

\item {Carleman, T}: {Les Fonctions Quasi-analytiques}.
 Gauthier-Villars, Paris (1926)

\item Chen, Y, Karagiannidis, GK, Lu, H, Cao, N: Novel approximations to the statistics of products of independent random variables
and their applications in wireless communications. IEEE Trans.
Veh. Tech. {\bf 61}, 443--454 (2012)

\item  {Chihara, TS:} On indeterminate
Hamburger moment problems. {Pacific J. Math.} {\bf 27}, 475--484
(1968)

\item  {Chihara, TS:} A characterization and a class of distribution functions for the Stieltjes--Wigert polynomials.
Canad. Math. Bull. {\bf 13}, 529--532 (1970)

\item {Chow, YS, Teicher, H}: {Probability Theory: Independence, Interchangeability, Martingales}. 3rd
edn. Springer, New York (1997)

\item Christiansen, JS: The moment problem associated with the Stieltjes--Wigert polynomials. J. Math. Anal. Appl. {\bf 277}, 218--245 (2003)

\item Fischer, H: A History of the Central Limit Theorems: From Classical to Modern Probability Theory. Springer, New York (2011)

\item {Fr$\acute{\hbox{e}}$chet, M, Shohat, J}: A
proof of the generalized second limit theorem in the theory of
probability. {Trans. Amer. Math. Soc.} {\bf 33}, 533--543 (1931)

\item {Galambos, J, Simonelli, I}: {Products of Random
Variables: Applications to Problems of Physics and to
Arithmetical Functions}. Marcel Dekker, New York (2004)

\item Garnett, JB: {Bounded Analytic Functions.}
Springer, New York (1981)

\item  Goffman, C, Pedrick, G:
 {First Course in Functional Analysis.} Prentice Hall of India,
 New Delhi (2002)

\item G\'omez, R, L\'opez-Garc\'ia, M: A family of heat
functions as solutions of indeterminate moment problems. {Int.
J. Math. Math. Sci.} Article ID 41526, 1--11 (2007)

\item Graffi, S, Grecchi, V: Borel summability and
indeterminacy of the Stieltjes moment problem: Application to
the anharmonic oscillators. {J. Math. Phys.} {\bf 19},
1002--1006 (1978)

\item Gut, A: On the moment problem. {Bernoulli} {\bf
8}, 407--421 (2002)

\item Gut, A: On the moment problem for random sums.
{J. Appl. Probab.} {\bf 40}, 797--802 (2003)

\item {Hamburger, H}: \"Uber eine Erweiterung des
Stieltjesschen Momentenproblems I--III. {Math. Ann.} {\bf 81},
235--319 (1920), {\bf 82}, 120--164 (1921), {\bf 82}, 168--187
(1921)

\item  {Hardy,
GH}: On Stieltjes' `probl\`{e}me des moments'. {Messenger of
Math.} {\bf 46}, 175--182 (1917), {\bf 47}, 81--88 (1918)
[Collected Papers of G.H. Hardy, Vol. VII, pp.\,75--83, 84--91
(1979)  Oxford University Press, Oxford.]

\item {Heyde, CC}: On a property of the lognormal distribution.
{J. Roy. Statist. Soc. Ser. B}  {\bf 25}, 392--393 (1963)

\item Heyde, CC: Some remarks on the moment problem (I).
{Quart. J. Math. Oxford} {\bf 14}, 91--96 (1963)

\item {H\"orfelt, P}: The moment problem for some Wiener
functionals: corrections to previous proofs (with an appendix by
H. L. Pedersen). {J. Appl. Probab.} {\bf 42}, 851--860 (2005)

\item Kjeldsen, TH: The early history of the moment
problem. {Historia Math.} {\bf 20}, 19--44 (1993)

\item {Klebanov, LB, Mkrtchyan, ST}: Estimation of the
closeness of distributions in terms of identical moments. In:
 {Stability Problems for Stochastic Models}, (Proc.\,Fourth
 All-Union Sem., Palanga, 1979) (Russian), Zobotarev, VM,
Kalashnikov, VV (eds.), pp.\,64--72, Moscow (1980).
Translations:  {J. Soviet Math.} {\bf 32}, 54--60 (1986);
{Selected Translations in Mathematical Statistics and
Probability} {\bf 16}, 1--10 (Estimating the proximity of
distributions in terms of coinciding moments) (1985)

\item {Kleiber, C}: On moment indeterminacy of the Benini
income distribution. {Statist. Papers} {\bf 54}, 1121--1130
(2013)

\item {Kleiber, C}: The generalized lognormal distribution
and the Stieltjes moment problem. {J. Theor. Probab.} {\bf 27},
1167--1177 (2014)

\item {Koosis, P}: {The Logarithmic Integral I.} Cambridge
Univ.\,Press, Cambridge (1988)

\item Kopanov, P, Stoyanov, J: Lin's condition for
functions of random variables and moment determinacy of
probability distributions. C. R. Bulg. Acad. Sci. {\bf 70},
611--618 (2017)

\item {Krein, M}: On a problem of extrapolation of A.N.
Kolmogoroff. {Comptes Rendus (Doklady) l'Academie Sci l'URSS}
{XLVI} (8), 306--309 (1945) [Dokl. Akad. Nauk SSSR {\bf 46},
339--342 (1944)]

\item Lariviere, MA: A note on probability
distributions with increasing generalized failure rates. {Oper.
Res.} {\bf 54}, 602--604 (2006)

\item Leipnik, R: The lognormal distribution and strong
non-uniqueness of the moment problem. {Theory Probab. Appl.}
{\bf 26}, 850--852 (1982)

\item {Lin, GD}: On the moment problems. {Statist. Probab.
Lett.} {\bf 35}, 85--90 (1997) Erratum: ibid {\bf 50}, 205
(2000)

\item Lin, GD, Huang, JS: The cube of a logistic
distribution is indeterminate. {Austral. J. Statist.} {\bf 39},
247--252 (1997)

\item Lin, GD, Stoyanov, J: On the moment determinacy of
the distribution of compound geometric sums. {J. Appl. Probab.}
{\bf 39}, 545--554 (2002)

\item {Lin, GD, Stoyanov, J}: The logarithmic skew-normal
distributions are moment-indeterminate. {J. Appl. Probab.} {\bf
46}, 909--916 (2009)

\item {Lin, GD, Stoyanov, J}: Moment determinacy of powers and
products of nonnegative random variables. {J. Theoret. Probab.}
{\bf 28}, 1337--1353 (2015)

\item {Lin, GD, Stoyanov, J}: On the moment determinacy of
products of non-identically distributed random variables.
{Probab. Math. Statist.} {\bf 36}, 21--33 (2016)

\item L\'opez-Garc\'ia, M: Characterization of solutions to the
log-normal moment problem. Theory Probab. Appl. {\bf 55},
303--307 (2011)

\item {Ostrovska, S}: Constructing Stieltjes classes for
M-indeterminate absolutely continuous probability distributions.
{ALEA, Lat. Am. J. Probab. Math. Stat.} {\bf 11}, 253--258
(2014)

\item Ostrovska, S: On the powers of polynomial logistic
distributions. {Braz. J. Probab. Stat.} {\bf 30}, 676--690
(2016)

\item {Ostrovska, S, Stoyanov, J}: Stieltjes classes for
M-indeterminate powers of inverse Gaussian distributions.
{Statist. Probab. Lett.} {\bf 71}, 165--171 (2005)

\item {Pakes, AG}: Length biasing and laws equivalent to the log-normal.
{J. Math. Anal. Appl.} {\bf 197}, 825--854 (1996)

\item {Pakes, AG}: Remarks on converse Carleman and Krein
criteria for the classical moment problem. {J. Aust. Math. Soc.}
{\bf 71}, 81--104 (2001)

\item {Pakes, AG}: Structure of Stieltjes classes of moment-equivalent
probability laws. {J. Math. Anal. Appl.} {\bf 326}, 1268--1290
(2007)

\item Pakes, AG: {On generalized stable and related laws}. {J. Math. Anal.
Appl.} {\bf 411}, 201--222 (2014)

\item Pakes, AG, Hung, W-L,  Wu, J-W: Criteria for the unique
determination of probability distributions by moments. {Aust.
N.Z. J. Statist.} {\bf 43}, 101--111 (2001)

\item Pakes, AG, Khattree, R: Length-biasing, characterizations of laws and the
moment problem. Austral. J. Statist. {\bf 34}, 307--322 (1992)

\item {Pedersen, HL}: On Krein's theorem for indeterminacy of
the classical moment problem. {J. Approx. Theory} {\bf 95},
90--100 (1998)

\item Penson, KA, Blasiak, P, Duchamp, GHE, Horzela, A, Solomon,
AI: On certain non-unique solutions of the Stieltjes moment
problem. {Discrete Math. Theor. Comput. Sci.} {\bf 12}, 295--306
(2010)

\item Prohorov, YuV, Rozanov, YuA:
{Probability Theory}. Translated by K.\,Krickeberg and
H.\,Urmitzer. Springer, New York (1969)

\item Rao, CR, Shanbhag, DN,  Sapatinas, T, Rao, MB:
Some properties of extreme stable laws and related infinitely
divisible random variables. {J. Statist. Plann. Inference} {\bf
139}, 802--813 (2009)

\item Serfling, RJ: {Approximation Theorems of
Mathematical Statistics.} Wiley, New York (1980)

\item {Shohat, JA, Tamarkin, JD}: {The Problem of
Moments.}  Amer. Math. Soc., New York (1943)

\item {Slud, EV}: The moment problem for polynomial forms in
normal random variables. {Ann. Probab.} {\bf 21}, 2200--2214
(1993)

\item Springer, MD: The Algebra of Random Variables. Wiley, New York (1979)

\item {Stieltjes, TJ}: Recherches sur les fractions
continues. {Ann. Fac. Sci. Univ. Toulouse Math.} {\bf 8(J)},
1--122 (1894); {\bf 9(A)}, 1--47 (1895). Also in: {Stieltjes,
T.J. Oeuvres Completes}. Noordhoff, Gr\"oningen {\bf 2},
402--566 (1918)

\item Stirzaker, D: The Cambridge Dictionary of Probability and
Its Applications. Cambr. Univ. Press, Cambridge (2015)

\item Stoyanov, J: Inverse Gaussian distribution and the moment
problem. {J. Appl. Statist. Sci.} {\bf 9}, 61--71 (1999)

\item Stoyanov, J: Krein condition in probabilistic
moment problems. {Bernoulli} {\bf 6}, 939--949 (2000)

\item {Stoyanov, J}: Stieltjes classes for moment-indeterminate
probability distributions. {J.  Appl. Probab.} {\bf 41},
281--294 (2004)

\item Stoyanov, JM: {Counterexamples in Probability}. 3rd edn.
 Dover Publications, New York (2013) [First and second edns:
Chichester: Wiley, 1987 and 1997.]

\item  Stoyanov, J, Kopanov, P: Lin's condition and moment
determinacy of functions of random variables. Revised for
{Statist. Probab. Lett.} (2017)

\item {Stoyanov, J, Lin, GD}: Hardy's condition in the moment
problem for probability distributions. {Theory Probab. Appl.}
{\bf 57},  699--708 (2013)

\item {Stoyanov, J, Lin, GD, DasGupta, A}: Hamburger moment
problem for powers and products of random variables.  {J.
Statist. Plann. Inference} {\bf 154}, 166--177 (2014)

\item {Stoyanov, J, Tolmatz, L}: New Stieltjes classes
involving generalized gamma distributions. {Statist. Probab.
Lett.} {\bf 69}, 213--219 (2004)

\item {Stoyanov, J, Tolmatz, L}: Method for constructing Stieltjes
classes for M-indeterminate probability distributions. {Appl.
Math. Comput.} {\bf 165}, 669--685 (2005)

\item  {Targhetta, ML}: On a family of indeterminate
distributions. {J. Math. Anal. Appl.} {\bf 147}, 477--479 (1990)

\item {Wang, J}: Constructing Stieltjes classes for power-order
M-indeterminate distributions. {J. Appl. Probab. Statist.} {\bf
7}, 41--52 (2012)

\end{description}

\end{document}